\documentclass[a4paper,11pt]{article}

\usepackage{array}
\usepackage{theorem}
\usepackage{amsmath,amscd,amssymb} 
\usepackage{latexsym}
\usepackage{epsfig}
\theorembodyfont{\sl}

\newtheorem{lemma}{Lemma}[section]

\newtheorem{proposition}[lemma]{Proposition}
\newtheorem{theorem}[lemma]{Theorem}
\newtheorem{corollary}[lemma]{Corollary}

\newtheorem{remark}[lemma]{Remark}

\renewcommand{\AA}{\mathbb A}

\newcommand{\CC}{\mathbb C}

\newcommand{\NN}{\mathbb N}

\newcommand{\PP}{\mathbb P}
\newcommand{\QQ}{\mathbb Q}
\newcommand{\RR}{\mathbb R}

\newcommand{\ZZ}{\mathbb Z}

\newcommand{\cI}{\mathcal I}

\newcommand{\cO}{\mathcal O}

\newcommand{\cQ}{\mathcal Q}

\newcommand{\To}{\longrightarrow}

\newcommand{\Sum}{\sum\limits}

\renewcommand{\Tilde}{\widetilde}
\renewcommand{\Hat}{\widehat}

\newcommand{\Pic}{\mathop{\mathrm {Pic}}\nolimits}
\newcommand{\Sing}{\mathop{\mathrm {Sing}}\nolimits}

\newcommand{\hcf}{\mathop{\mathrm {hcf}}\nolimits}

\newcommand{\Area}{\mathop{\null\mathrm {Area}}\nolimits}

\newcommand{\gothm}{\mathfrak m}

\newcommand{\Pw}{\PP^4(w)}
\newcommand{\hX}{\Hat{X}}
\newcommand{\hd}{\Hat{d}}
\newcommand{\hs}{\Hat{s}}
\newcommand{\hchi}{\chi\left(\cO_{\hX}\right)}
\newcommand{\hpi}{\Hat{\pi}}

\newcommand{\hdelta}{\Hat{\delta}}

\newcommand{\hH}{\Hat{H}}
\newcommand{\hK}{{K_{\hX}}}
\newcommand{\tX}{{\Tilde{X}}}
\newcommand{\tchi}{\chi\left(\cO_{\tX}\right)}

\begin{document}

\pagenumbering{arabic}

\setcounter{section}{0}
\title{Boundedness for surfaces in weighted $\PP^4$}
\author{L.V. Rammea \& G.K. Sankaran}
\maketitle
\begin{abstract}
 Ellingsrud and Peskine (1989) proved that there exists a bound on the
 degree of smooth non general type surfaces in $\PP^4$. The latest
 proven bound is 52 by Decker and Schreyer in 2000.

 In this paper we consider bounds on the degree of a quasismooth
 non-general type surface in weighted projective $4$-space. We show
 that such a bound in terms of the weights exists, and compute an
 explicit bound in simple cases.
\end{abstract}

\section*{Introduction}

Ellingsrud and Peskine~\cite{EP} proved that there exists an integer
$d_{0}$ such that all smooth non-general type surfaces in $\PP^4$ have
degree less than or equal to $d_{0}$. This motivated a search for such
surfaces, partly by computational methods, and also an effort to find
an effective bound on $d_0$, begun by Braun and Fl\o{}ystad
in~\cite{BF}.  As far as we know the smallest proven bound is $52$ by
Decker and Schreyer \cite{DS}.

Some of the methods used to find such surfaces are also applicable to
surfaces in weighted projective spaces $\Pw$ (some first steps in this
direction are taken in~\cite{Ra}). It is therefore natural
to ask whether a similar bound can be found for the degree of
quasismooth non-general type surfaces in a weighted projective space
with given weights. In this paper we show that such a readily
computable bound (of course depending on the weights) does exist, and
we compute it in some cases.

To show that a bound exists all we need is a fairly simple adaptation
of the way in which the results of~\cite{EP} (or~\cite{BF}) are
applied. For a computable bound we use the results of~\cite{BF}
together with some information about the contribution from the
singularities of the surface in~$\Pw$.

Our procedure is to exploit the representation of $\Pw$ as a quotient
of $\PP^4$ by a finite group action. Starting with a quasismooth
non-general type surface $X$ in weighted projective $4$-space $\Pw$,
we take its cover in $\PP^4$. This will (usually) be of general type,
but it will have invariants bounded in terms of those of $X$, and the
results of~\cite{BF} still apply in this situation. 

\emph{Acknowledgements.} This work forms part of the Bath Ph.D. thesis
\cite{Ra} of the first author, supported by a Commonwealth Scholarship
of the Association of Commonwealth Universities.

\section{Bounding the degrees}

We fix weights $w=(w_0,w_1,w_2,w_3,w_4)$ with $w_i\in\NN$: unless
otherwise stated, $i$ and $j$ always denote indices in the range $0\le
i,\ j\le 4$. We may assume that any four of the $w_i$ are coprime:
such weights are called well-formed (see \cite[1.3.1]{Do} and
\cite[5.9 \& 5.11]{IF}).
We also order the weights so that $w_i\le w_{i+1}$: in particular, the
largest weight is~$w_4$.  We write $|w|$ for the sum of the weights,
and $m$ for their product.
%
%
The weighted projective space $\Pw$ of dimension $4$ is defined to be
the quotient $(\CC^5\setminus\{0\})/\CC^*$, where $\CC^*$ acts by
\[
t\colon (x_0,\ldots,x_4)\to (t^{w_0}x_0,\ldots,t^{w_4}x_4).
\]
A surface $X\subset \Pw$ is said to be quasismooth if its punctured
affine cone $X^*$ is smooth: that is, if $X^*=q^{-1}(X)$ is smooth,
where $q\colon(\CC^5\setminus\{0\})/\CC^*\to \Pw$ is the quotient map:
see~\cite[3.1.5]{Do} or \cite[6.3]{IF}.

Alternatively (\cite[1.2.2]{Do}) we may regard $\Pw$ as a quotient of
$\PP^4$ under an action of the group $G_w=\prod_i \ZZ/w_i\ZZ$ of order
$m$. A generator $g_i$ of the $i$th factor acts by $x_i\mapsto
x_i^{w_i}$. We denote the quotient map $\PP^4\to \Pw$ by $\phi_w$.

Suppose that $X$ is a quasismooth surface, not of general type, in
$\Pw$.  Denote by $\hX$ the cover of $X$ in $\PP^4$ under the
$m$-to-$1$ map $\phi_w$: then $\hX$ is smooth. We always assume
that $X$ and $\hX$ are nondegenerate: that is, $\hX$ is not contained
in any hyperplane in $\PP^4$.

Let $f\colon \Tilde X\to X$ be the minimal resolution of $X$ (note
that $\Tilde X$ need not be a minimal surface).
\[
\begin{CD}
{} @. \hX @.\ \subset\  @. \PP^4\\
@. @V{\phi_w}VV @. @VV{\phi_w}V\\
\tX @>{f}>> X @.\ \subset\  @. \Pw
\end{CD}
\]
Further let $d$ be the degree of $X\subset \Pw$ and $\pi$ the
sectional genus of $X$. These are defined as follows: $\Pw$ and $X$
are $\QQ$-factorial varieties and there are $\QQ$-line bundles
$\cO_{\Pw}(1)$, $\cO_X(1)$ and $K_X$. Writing $H$ for the class of
$\cO_X(1)$ in $\Pic X\otimes \QQ$ and using the intersection form on
$\Pic X$ we have $d=H^2$ and $2\pi-2=H\cdot(H+K_X)$, so $d,\
\pi\in\QQ$.

We let $\hd$ be the degree of $\hX$ and $\hpi$ the sectional genus of
$\hX$. We put
\[
\hs=\min\left\{k|h^{0}\cI_{\hX}(k)\neq 0 \right\}
\]
and denote by $\sigma_f$ the number of irreducible exceptional curves of
$f$.

We first collect the facts about these invariants of the smooth
surface $\hX\subset\PP^4$.

\begin{proposition}\label{basic_bounds}
If $\hX\subset \PP^4$ is a smooth surface (possibly of general
type), and $r\le \hs$ and $r^2<\hd$, then
\begin{equation}\label{pihat_upper_bound}
2\hpi \leq \dfrac{\hd^2}{r}+(r-4)\hd+1.
\end{equation}
Moreover
\begin{equation}\label{doublepoint} 
\hd^2 -5 \hd -10(\hpi-1) +  12\hchi-2K^2_{\hX}=0.
\end{equation}
Finally, if $\hd>\hs(\hs-1)$ we have the lower bound for $\hchi$ 
\begin{eqnarray}\label{chi_lower_bound}
\hchi&\ge& \frac{\hd^3}{6\hs}+\hd^2\frac{\hs-5}{4\hs}+
\hd\frac{3\hs^2-30\hs+71}{24}\\
&&-\frac{\hs^4-5\hs^3-\hs^2+5\hs}{24}
-\frac{\gamma^2}{2}-\gamma(\frac{\hd}{\hs}+s-\frac{5}{2})\nonumber
\end{eqnarray}
where $0\le \gamma\le \hd(\hs-1)^2/2\hs$.
\end{proposition}
\emph{Proof:} The inequality \eqref{pihat_upper_bound} is a
consequence of \cite[(B),\ (C), page 2]{EP}.  Let $\hH$ denote a
general hyperplane section of $\hX$, so that
$\hpi=g(\hH)$. According to \cite{Ro} (as quoted in \cite[(C),
page 2]{EP}), if $\hs>r$ and $\hd>r^2$ then $\hH\subset \PP^3$
does not lie on any surface of degree $<r$. Therefore, according to
\cite[(B), page 2]{EP}, we have $r(2\hpi-2)\le \hd^2+r(r-4)\hd$.  If
$\hs=r$ then (again by \cite[(B), page 2]{EP}) we have the same
inequality because then $\hH$ does lie on a surface of degree~$r$.

Equation~\eqref{doublepoint} is the double point formula as stated
in~\cite{EP} and \cite{BF}.
The estimate
\eqref{chi_lower_bound} is~\cite[(1.1)(e)]{BF}.  \smallskip

A more precise version of~\eqref{pihat_upper_bound}, valid under
certain conditions, is given in~\cite[(1.1)]{BF}.  In order to bound
the degree of smooth surfaces in $\PP^4$ what is needed is not the
precise form of \eqref{chi_lower_bound} but an estimate of the form
$\hchi\ge a(\hs)\hd^3+o(\hd^3)$, where $a(\hs)$ is some positive
constant depending on $\hs$ only. Ellingsrud and Peskine proved the
existence of such a bound in~\cite{EP} but did not give an explicit
one.

It will be convenient to work with the invariants $c_1^2(S)=K_S^2$ and
$c_2(S)$ (which is the topological Euler number $e(S)$) of a smooth
projective surface~$S$: these are connected by Noether's formula
\begin{equation}\label{noether}
12\chi(\cO_S)=c_1^2(S)+c_2(S)
\end{equation}
Since we are assuming that $\tX$ is not of general type we have (as in
\cite{EP} and \cite{BF}) that $K^2_{\tX}\leq 9$. Moreover, unless
$\tX$ is a rational surface with $K^2_{\tX}\ge 6$ we also have $6\tchi
\geq K^2_{\tX}$ (i.e.\ $c_2(\tX)- c_1^2(\tX)\ge 0$). If $\tX$ is a
rational surface then $\tchi=1$ so $c_2(\tX)-c_1^2(\tX) =
12\tchi-2K^2_{\tX} = 12-2K_{\tX}^2\ge -6$. So in any case if $X$ is
not of general type we have
\begin{equation}\label{c1squaredc2}
c_1^2(\tX)-c_2(\tX)\le 6.
\end{equation}
So we need to estimate $\hd$ and $\hpi$ in terms of $d$ and $\pi$, and
$K^2_{\hX}$ and $\hchi$ in terms of $K^2_{\tX}$ and $\tchi$. We shall
show the two propositions below.
\begin{proposition}\label{boundchernnumbers}
Suppose $X$ is a quasismooth normal surface in $\Pw$. Then
\begin{equation}\label{theta1}
c_1^2(\hX)\le mc_1^2(\tX)+\theta_1
\end{equation}
where 
\begin{equation}\label{theta1form}
\theta_1=k_0+k_1\hd+k_2\hdelta
\end{equation}
for suitable $k_0,\ k_1,\ k_2$ depending only on the weights
$w_i$. Moreover
\begin{equation}\label{theta2}
c_2(\hX)\ge mc_2(\tX)-\theta_2,
\end{equation}
and
\begin{equation}\label{sumtheta}
\theta_1+\theta_2=k'_0+k'_1\hd+k'_2\hdelta
\end{equation}
for suitable $k'_0,\ k'_1,\ k'_2$ depending only on the weights, and
$k_2'>-5$.
\end{proposition}

This proposition will be proved in Sections~\ref{boundc1} and
\ref{boundc2}, below.

Our main qualitative result is then the following.
\begin{theorem}\label{boundexists}
  There exists $d_w\in\NN$ depending only on the weights $w_i$ such
  that any quasi-smooth normal surface $X\in\Pw$ of degree $d>d_w$ is
  of general type.
\end{theorem}
\emph{Proof:}
We have seen that $\hX \to X$ is $m$-to-$1$, so 
\begin{equation} \label{hatd_is_md} 
\hd = md 
\end{equation}
so it is sufficient to show that if $X$ is not of general type then
$\hd$ is bounded by a function of the weights.

Suppose then that $X$ is not of general type.  We have, by adjunction,
$2\hpi-2=\hH\cdot(\hH+\hK)=\hd+\hdelta$, where $\hH$ is a hyperplane
section of $\hX$. Therefore by the estimate \eqref{pihat_upper_bound}
we obtain
\begin{equation}\label{deltahat}
\hdelta\le\frac{1}{r}\hd^2+(r-5)\hd
\end{equation}
as long as $r\le\hs$ and $r^2<\hd$.  We may also write the double
point formula as
\begin{equation}\label{doublepoint2}
\hd^2-10\hd-5\hdelta+c_2(\hX)-c_1^2(\hX)=0.
\end{equation}
By Proposition~\ref{boundchernnumbers} and the
inequality~\eqref{c1squaredc2} we have
\begin{equation}\label{c1squaredc2hat}
c_2(\hX)-c_1^2(\hX)\ge -6m-(\theta_1+\theta_2),
\end{equation}
so
\begin{eqnarray}\label{dhatsquaredbound}
0&\ge& \hd^2-10\hd-5\hdelta-6m-(\theta_1+\theta_2)\nonumber\\
 &=&\hd^2-(10+k'_1)\hd-(6m+k'_0)-(5+k'_2)\hdelta.
\end{eqnarray}

Combining this with \eqref{deltahat} gives (since $5+k'_2>0$)
\begin{eqnarray*}
0&\ge&\hd^2-(10+k'_1)\hd-(6m+k'_0)-(5+k'_2)(\frac{1}{r}\hd^2+(r-5)\hd)\\ 
&=&\left(1-\frac{5+k'_2}{r}\right)\hd^2
-((10+k'_1+(5+k'_2)(r-5))\hd-(6m+k'_0).
\end{eqnarray*}
So if $\hs>k'_2+5$ we may take $r=k'_2+6$ and this bounds $\hd$ in
that case.

On the other hand, suppose that $X$ is not of general type and $\hs\le
k'_2+5$. Then using Noether's formula, the double point formula
\eqref{doublepoint2}, and \eqref{chi_lower_bound} we have
\begin{eqnarray*}
0&=& \hd^2-10\hd-5\hdelta+12\hchi-2c_1^2(\hX)\\
&\ge& -2c_1^2(\hX)+ \frac{2}{\hs}\hd^3+O(\hd^2)\\
&\ge& -2mc_1^2(\tX)-\theta_1+\frac{2}{\hs}\hd^3+O(\hd^2)\\
&\ge& \frac{2}{\hs}\hd^3+O(\hd^2)-18m-k_0-k_1\hd-k_2\hdelta\\
&=& \frac{2}{\hs}\hd^3+O(\hd^2)
\end{eqnarray*}
by \eqref{theta1form} and \eqref{deltahat}: the constants depend on
$\hs$ but this is now bounded in terms of the weights. So again we
obtain a bound for $\hd$ in terms of the $w_i$.

\section{Singularities of $\Pw$ and of $X$}

In this section we collect some preliminary information about the
action of $G_w$ on $\PP^4$ and on $\hX$. We choose an isomorphism
$G_w\to \prod\ZZ/w_i\ZZ$ by choosing generators $g_i\in G_w$ of order
$w_i$. The singularities arise at fixed points of the $G_w$-action, so
let us consider those.

Suppose that $x=(x_0:\ldots:x_4)\in\PP^4$ is fixed by
$g=g_0^{a_0}\ldots g_4^{a_4}$. Without loss of generality we take
$x_0=1$: then for $j\neq 0$ we have $\zeta_0^{-a_0}\zeta_j^{a_j}=1$,
where $\zeta_j=e^{2\pi i/w_j}$.

\begin{lemma}\label{fixed_points} 
  If $x\in\PP^4$ is fixed by a non-trivial element of $G_w$, then $x$
  lies in a coordinate linear subspace $\PP_J$ given by
  $\PP_J=\{x_j=0\mid j\in J\subset\{0,\ldots,4\}\}$. The stabiliser
  of a general point of $\PP_J$ is the group $\Gamma_J$ generated by the
  $g_j$ for $j\in J$ and the element $g_J=\prod_{i\not\in
    J}g_i^{w_i/r_J}$, where $r_J=\hcf(a_i\mid i\not\in J)$.
\end{lemma}

This is immediate from the description of the action above. By a
general point in $\PP_J$ is meant, in this case, a point that is not
in $\PP_{J'}$ for any $J'\supset J$.

\begin{lemma}\label{Xsings}
  The singularities of $X$ are cyclic quotient singularities whose
  order divides one of the weights.
\end{lemma}
\emph{Proof:} At a fixed point $x\in\PP^4$, the elements $g_j\in
\Gamma_J$ act on the tangent space by quasi-reflections: the $j$th
eigenvalue is $\zeta_j^{a_j}$ and the others are~$1$. So the quotient
by the subgroup $\Gamma'_J$ generated by those elements is smooth, and the
singularity of $\PP_w$ or of $X$ at $z=\phi_w(x)$ is a quotient by the
action of the cyclic group generated by $g_J$. The order of this
element, or of its image in $\Gamma_J/\Gamma'_J$, is $r_J$, which
divides $w_i$ for $i\not\in J$.
 
\begin{remark}\label{smooth_hyperplanes}
  If $\#J=1$ then $r_J=1$ since the weights are well-formed, so the
  general point of a coordinate hyperplane in $\Pw$ is smooth. For
  each $i$, the number of singular points of $X$ with $z_i=0$ is at
  most $\hd$.
\end{remark}
\begin{remark}\label{coprime_singularities}
  If the weights are pairwise coprime then the singularities occur at
  the points $P_0=(1:0:\ldots:0),\ldots,P_4=(0:\ldots:0:1)\in\Pw$, and
  the singularity of $\PP_w$ at $P_i$ has order exactly $w_i$. If $X\ni
  P_i$ then $X$ also has a cyclic quotient singularity of order $w_i$
  at $P_i$.
\end{remark}

\begin{lemma}\label{transversality}
Suppose that $(Y,0)$ is a nondegenerate smooth surface germ in
$(\AA^4,0)$ with coordinates $t_1,\ldots,t_4$ at $0\in\AA^4$. Let
$\gamma$ be the quasi-reflection $\gamma(t_1)=\xi t_1$, where $\xi$ is
a primitive $n$th root of unity, and that $Y$ is
$\gamma$-invariant. Then $Y$ meets $A=(t_1=0)$ tranversely.
\end{lemma}
\emph{Proof:} Suppose not: then $T_{Y,0}\subset A$. Therefore the
ideal $\cI_{Y,0}\subset \cO_{\AA^4,0}$ contains an element $f$ of the
form $f=t_1+h$ with $h\in \gothm^2\subset \cO_{\AA^4,0}$, where
$\gothm$ is the maximal ideal of $\cO_{\AA^4,0}$.

We write $h=\Sum_{\nu=0}^{n-1}h_\nu$, where
$\gamma(h_\nu)=\xi^\nu(h_\nu)$: if we write $h$ as a polynomial in
$t_1$, so  $h=\Sum_ra_r(t_2,t_3,t_4)t_1^r$, we have
$h_\nu=\Sum_{r\equiv \nu \mod n}a_r(t_2,t_3,t_4)t_1^r$. Then
$$
\cI_{Y,0}\ni (1+\gamma+\gamma^2+\cdots+\gamma^{n-1})(f)=nh_0
$$
so $\cI_{Y,0}\ni f-h_0=t_1+\Sum_{\nu\neq 0}h_\nu$. But
$t_1$ divides the right-hand side, so since $h\in \gothm^2$ we have
$f-h_0=t_1(1+b)$, where $b\in \gothm$. Since $\cI_{Y,0}$ is a
prime ideal contained in $\gothm$ this implies $t_1\in \cI_{Y,0}$,
contradicting the nondegeneracy.  
\begin{corollary}\label{smoothbranchcurves}
  If $w_i\neq 1$, then $\hX$ meets the ramification divisor
  $\PP_{\{i\}}$ transversely and the curve ${\Hat C}_i=\hX\cap
  \PP_{\{i\}}$ is a smooth curve of genus $\hpi$.
\end{corollary}
\emph{Proof:} The second part follows immediately from the
first, which is immediate from Lemma~\ref{transversality}.

\section{Comparing $c_1^2$.}\label{boundc1}

In this section we prove \eqref{theta1} and \eqref{theta1form} from
Proposition~\ref{boundchernnumbers}, and give values for the constants
$k_0$, $k_1$ and $k_2$.

Let $\Delta=\sum_{1\le \nu\le \sigma_f} a_\nu E_\nu$ be the
discrepancy of $f$, so that $a_\nu\in\QQ$ and $K_{\tX}=f^*K_X+\Delta$.
Then $f^*K_X\cdot \Delta$ vanishes and $(f^*K_X)^2=K_X^2$, so
$K_{\tX}^2= K_X^2+\Delta^2$.

\begin{lemma}\label{bound_Delta2}
  If $f_0\colon \Tilde Y\to Y$ is the minimal resolution of a isolated
  cyclic quotient $(Y,0)$ of order $n$ and the discrepancy of $f_0$ is
  $\Delta_0$, then $0> \Delta_0^2\ge -n$.
\end{lemma}
\emph{Proof:}
 This (which is not a sharp bound) is most easily seen by
toric methods. If the singularity is $\frac{1}{n}(1,a)$ with $(n,a)=1$
then the minimal resolution is described by taking the decomposition
given by the convex hull of $\ZZ^2+\frac{1}{n}(1,a)\ZZ$ in the first
quadrant of $\RR^2$. The exceptional curves $E_\nu$, $0<\nu<k$, correspond to
primitive vectors $P_\nu=(x_\nu,y_\nu)$ of this lattice: put
$\ell_\nu=x_\nu+y_\nu$, and write $E_0$ and $E_k$ for the toric curves
corresponding to the rays spanned by $(1,0)$ and $(0,1)$. Then we have
$E_\nu E_{\nu\pm 1}=1$ and $E_\mu E_\nu=0$ if $\mu\neq \nu,\ \nu\pm
1$. Moreover on $\tilde Y$ we have $\Sum_{0\le \nu\le k}\ell_\nu
E_\nu\equiv 0$ (linear equivalence), and
$\Delta=-\Sum_{0<\nu<k}E_\nu$. Therefore 
\begin{eqnarray*}
\Delta_0^2&=&\sum_{0<\nu<k}E_\nu(\sum_{0<\mu<k}E_\mu)\\
&=&\sum_{0<\nu<k}E_\nu\Big((\sum_{\mu\neq
  0,\ \nu,\ k}E_\mu)+E_\nu\Big)\\
&=&\sum_{0<\nu<k}E_\nu(-E_0-E_k+\sum_{\mu\neq \nu}
(1-\frac{\ell_\mu}{\ell_\nu})E_\mu)\\
&=&-2-\sum_{0<\nu<k}((\frac{\ell_{\nu-1}}{\ell_\nu}-1)
+(\frac{\ell_{\nu+1}}{\ell_\nu}-1))
\end{eqnarray*}

Suppose for definiteness that $\ell_{nu+1}>\ell_{\nu}$. Then
$\frac{\ell_{\nu+1}}{\ell_\nu}-1$ is twice the area (relative to the
lattice $\Lambda=\ZZ^2+\frac{1}{n}(1,a)\ZZ$) of the triangle
$T^+_\nu=P_\nu Q_\nu P_{\nu+1}$, where
$Q_{\nu}=\frac{\ell_{\nu+1}}{\ell_\nu}P_\nu$, since $\Area(OP_\nu
P_{\nu+1})=\frac{1}{2}$ relative to $\Lambda$. So
\begin{eqnarray*}
-\frac{1}{2}\Delta_0^2&\le& -1-\sum_{0<\nu<k}\Area(T_\nu^+)\\
&=&-\Area(OP_0P_1)-\Area(OP_{k-1}P_k)-\sum_{0<\nu<k}\Area(T_\nu^+).
\end{eqnarray*}
But these triangles do not overlap and they are contained in the unit
triangle $OP_0P_k$, which has area $\frac{n}{2}$ relative to $\Lambda$.

Now we compute $K_{\hX}^2$ from $K_{\hX}=\phi_w^{*}(K_X)+\sum(w_i
-1)\hH_i$, where $\hH_i=\PP_{\{i\}}\cap \hX=(x_i=0)$ and so we get

\begin{equation}\label{K2hatX_K2X}
K^2_{\hX}= m K^2_X+2\sum(w_i-1)\hdelta -\sum(w_i-1)(w_j-1)\hd
\end{equation}
since $\phi_w^*(K_X)^2=mK_X^2$.

\begin{proposition}\label{have_theta1}
We have $c_1^2(\hX)\le mc_1^2(\tX)+\theta_1$, where (recall that $w_4$
is the largest weight) 
\begin{equation}\label{theta1_formula}
\theta_1=\Big(10mw_4-\sum_{0\le i,j\le
  4}(w_i-1)(w_j-1)\Big)\hd + 2(|w|-5)\hdelta. 
\end{equation}
\end{proposition}
\emph{Proof:} For a singular point $z\in \Sing(X)$ we denote the
discrepancy at $z$ by $\Delta_z$. If $z\in H_J=\phi_w(\PP_J)\cap X$
then the order of the singularity is $r_J=\hcf(w_i\mid i\not\in
J)$. There are at most $\binom{5}{2}\hd$ distinct points on
the $\hH_{\{ij\}}$ altogether, so the total number of singular points is
at most $10\hd$.

Each singular point has order $r_J$ dividing some of the $w_i$, so
$\Delta_z^2\ge -r_J \ge -w_4$. Then
\begin{eqnarray*}
c_1^2(\tX)=K_{\tX}^2 &=& K_X^2+\Delta^2\\
&=& K_X^2+\sum_{z\in \Sing(X)}\Delta_x^2\\
&\ge& K_X^2-10w_4\hd.
\end{eqnarray*}
Now, using~\eqref{K2hatX_K2X}, we get
\begin{eqnarray*}
c_1^2(\hX)&=& mK_X^2 +2\hdelta\big(|w|-5\big)
-\hd\sum_{0\le i,j\le 4}(w_i-1)(w_j-1)\\
&\le& mc_1^2(\tX)+2\hdelta\big(|w|-5\big)
+\hd\Big(10mw_4-\sum_{0\le i,j\le 4}(w_i-1)(w_j-1)\Big)
\end{eqnarray*}
as required.

If the $w_i$ are pairwise coprime we can do slightly better. In that
case the only singularities are at the points $P_i$ if they are in
$X$. Therefore we have
\begin{equation}\label{K2isopoints}
c_1^2(\tX)=K_X^2+\sum_{P_i\in X}\Delta_i^2\ge K_X^2-\sum_i q_i w_i,
\end{equation}
where $\Delta_i$ is the discrepancy at $P_i$ and $q_i=1$ if $P_i\in
X$, $q_i=0$ if $P_i\not\in X$. This gives
\begin{equation}\label{K2hatisopoints}
c_1^2(\hX) \leq m c_1^2(\tX) + m \sum q_iw_i + 2
(|w|-5)\hdelta -\hd\sum_{0\leq i,j\leq 4} (w_i-1)(w_j-1).
\end{equation} 

\section{Comparing $c_2$}\label{boundc2}

Recall that if $x\in \hX\cap \PP_J$ then $\Gamma_J$ stabilises $x$. We put
$\hX_J =\hX\cap (\PP_J\setminus\bigcup_{J'\supset J}\PP_{J'})$. On
$\hX_J$ the stabiliser is precisely $\Gamma_J$. The
order of $\Gamma_J$ is $h_J=r_J\prod_{j\in J}w_j$: in particular,
$h_\emptyset=1$ and $h_{\{i\}}=w_i$.

$\hX_{\{i\}}$ is the complement of up to $4\hd$ points on a smooth
curve of genus $\hpi$, by Corollary~\ref{smoothbranchcurves}. Those
points lie in some $\hX_J$ with $\# J\ge 2$: in particular they all
lie on $\hH_j$ for some $j\neq i$, and there are $\hd$ such points for
each such $j$. They may not all be distinct, however. Therefore
\begin{equation}\label{eHi}
2-2\hpi>e(\hX_{\{i\}})\ge 2-2\hpi-4\hd.
\end{equation}
Denote by $\cQ$ the set of points of $\hX$ lying in at least two
coordinate hyperplanes of $\PP^4$: thus $\cQ=\hX\cap \bigcup_{\#
  J\ge 2}\PP_J$ as a set. The set $\cQ$ is finite, of cardinality $q\le
10\hd$, and $\hX=\hX_\emptyset\coprod \bigcup_i \hX_{\{i\}}\coprod
\cQ$.

We put $X_J=\phi_w(\hX_J)$, for $J\subset \{0,\ldots,4\}$, so that
\begin{equation*}
\phi_w|_{\hX_J}\colon \hX_J\To X_J
\end{equation*}
is unramified and its degree is $\lvert G_w: \Gamma_J\rvert=m/r_J$.

\begin{lemma}\label{c2tilde_c2}
For each $x\in \cQ$, let $r_x$ be the order of the singularity of
$z=\phi_w(x)\in X$, so $r_x=r_J$ if $x\in \hX_J$. Then 
\begin{equation*}
c_2(\tX)\le e(X)+\sum_{x\in\cQ} (r_x-1).
\end{equation*}
\end{lemma}
\emph{Proof:} The resolution $f\colon \tX\to X$, in a neighbourhood of
$z$, consists of a sequence of at most $r_x-1$ blow-ups, needed to
resolve the quotient singularity of order $r_x$ at $z\in X$. Therefore
$\sigma_f\le \sum_{x\in\cQ}(r_x-1)$. Each blow-up contracts a smooth
rational curve: topologically, therefore, $f$ contracts $\sigma_f$
$2$-spheres to points, and each of these contractions reduces the Euler
characteristic by~$1$, so $e(\tX)=e(X)+\sigma_f\le e(X)+\sum_{x\in
  \cQ}(r_x-1)$.

\begin{proposition}\label{have_theta2}
We have $c_2(\hX)\ge mc_2(\tX)-\theta_2$, where (recall that $w_4$
is the largest weight) 
\begin{equation}\label{theta2_formula}
\theta_2=\Big(10mw_4-(|w|-5)\Big)\hd - (|w|-5)\hdelta. 
\end{equation}
\end{proposition}
\emph{Proof:} By the additiviity of Euler characteristic we have
\begin{eqnarray*}
c_2(\hX)&=& \sum_J e(\hX_J)\\
&=& \sum_J \lvert G_w:\Gamma_J\rvert e(X_J)\\
&=& me(X_\emptyset)+\sum_{J\neq \emptyset}\lvert G_w:\Gamma_J\rvert e(X_J)\\
&=& m\Big(e(X)-\sum_{J\neq \emptyset}e(X_J)\Big)+\sum_{J\neq
  \emptyset}\lvert G_w:\Gamma_J\rvert e(X_J)\\
&=& me(X)+\sum_{J\neq \emptyset}(1-h_J)\lvert G_w:\Gamma_J\rvert e(X_J)\\
&=& me(X)+\sum_{J\neq \emptyset}(1-h_J) e(\hX_J).
\end{eqnarray*}
Write $h_x=h_J$ if $x\in \hX_J$. Using $\hX=\hX_\emptyset\coprod \bigcup_i \hX_{\{i\}}\coprod
\cQ$ and Lemma~\ref{c2tilde_c2}, this gives
\begin{eqnarray*}
c_2(\hX)&=& me(X)-\sum_i (w_i-1)e(\hX_{\{i\}})-\sum_{x\in\cQ}(h_x-1)\\
&\ge& mc_2(\tX)-\sum_i (w_i-1)e(\hX_{\{i\}})-m\sum_{x\in\cQ}(r_x-1)-\sum_{x\in\cQ}(h_x-1)\\
&\ge& mc_2(\tX)-(|w|-5)(2-2\hpi)-m\sum_{x\in\cQ}(r_x-1)-\sum_{x\in\cQ}(h_x-1)\\
&=& mc_2(\tX)+(|w|-5)(\hd+\hdelta)-(m+1)\sum_{x\in\cQ}(r_x-1)\\
&\ge& mc_2(\tX)+(|w|-5)(\hd+\hdelta)-10mw_4\hd,
\end{eqnarray*}
as claimed, since $q\le 10\hd$, $r_x\le w_4$ and $h_x\le m$.

We can now complete the proof of Proposition~\ref{boundchernnumbers}
and hence of Theorem~\ref{boundexists}, by remarking that from
Propositions~\ref{have_theta1} and~\ref{have_theta2} we get
\begin{equation*}
\theta_1+\theta_2=\Big(20mw_4-(|w|-5)-\sum(w_i-1)(w_j-1)\Big)\hd
+(|w|-5)\hdelta
\end{equation*}
so $k'_2=|w|-5>-5$.

\section{Examples}

It would of course be possible to obtain an explicit bound as in
Theorem~\ref{boundexists} from the argument above. However, such a
bound would be likely to be rather poor. In specific cases it is
possible to obtain a bound better than the general one implied
above. Although we still do not expect such a bound to be good, in the
sense that we expect that in fact all non general type surfaces will
be of much lower degree, in some cases it is not absurdly big.%
\medskip

\noindent\emph{Example 1: weights $(1,1,1,1,2)$}

We calculate a bound for the case of weights $(1,1,1,1,2)$. In this
case there is at most one singular point of $X$ and if there is a
singular point it is an ordinary double point. We let $q$ be the
number of singularities of $X$, so $q=0$ or $q=1$.

In this case the singularity, if any, is canonical and blowing up once
gives a crepant resolution, so $\Delta^2=0$ and
$c_1^2(\tX)=K_X^2$. Moreover $K_{\hX}=\phi^*K_X+\hH$, so 
\begin{eqnarray*}
c_1^2(\hX)&=&(\phi^* K_X+\hH)^2\\
&=&2K_X^2+2\phi^*K_X\hH+\hH^2\\
&=&2c_1^2(\tX)+2(K_{\hX}-\hH)\hH+\hH^2\\
&=&2c_1^2(\tX)-\hd+2\hdelta.
\end{eqnarray*}
We also have $c_2(\tX)=e(X)+q$ and 
\begin{eqnarray*}
c_2(\hX)&=&2e(X)-\sum_i(w_i-1)e(\hX_{\{i\}})-\sum_{x\in\cQ}(h_x-1)\\
&=&2e(X)-(2-2\hpi)-q\\
&=&2c_2(\tX)+\hd+\hdelta-3q.
\end{eqnarray*}
Thus $\theta_1=-\hd+2\hdelta$ and $\theta_2=3q-\hd-\hdelta$. Therefore
$k_0'=3q$, $k_1'=-2$ and $k_2'=1$, and \eqref{dhatsquaredbound} and
the formula below it give
$$
0\ge \left(1-\frac{6}{r}\right)\hd^2-(6r-22)\hd-(12+3q)
$$
as long as $r\ge\hs\ge 7$ and $r^2<\hd$. Taking $r=7$, we see that
$\hd\le 140$ in this case. (By taking $r=9$ we can obtain $\hd\le 96$,
but as we shall see that will not yield a better bound in the
end. Clearly taking $r\ge 10$ we cannot do better than $d\le 100$
because for this case we need $r^2<\hd$.)

We must also deal with the cases $\hs<7$: if we use $r=9$ we must also
handle $\hs=7$ and $\hs=8$ separately. But now we have, using
$c_1^2(\tX)\le 9$, the estimate \eqref{pihat_upper_bound} for
$\hdelta$ with $r=\hs$, the bounds on $\hchi$ and $\gamma$ from
Proposition~\ref{basic_bounds}, and the double point formula
\begin{eqnarray*}
0&=&\hd^2-10\hd+12\hchi-2c_1^2(\hX)\\
&=&\hd^2-10\hd+12\hchi-4c_1^2(\tX)+2\hd=4\hdelta\\
&\ge&\hd^2-8\hd+12\hchi-36-\frac{4}{\hs}\hd^2-4(\hs-5)\hd+4\\
&\ge&12\Big[\frac{\hd^3}{6\hs}+\hd^2\Big(\frac{\hs-5}{4\hs}\Big)
+\hd\Big(\frac{3\hs^2-30\hs+71}{24}\Big)
-\frac{\hs^4-5\hs^3-\hs^2+5\hs}{24}\\
&&-\frac{1}{2}\Big(\frac{(\hs-1)^4}{4\hs^2}\Big)\hd^2
-\Big(\frac{(\hs-1)^2}{2\hs^2}\Big)\hd^2
-\Big(\frac{(\hs-5/2)(\hs-1)^2}{2\hs}\Big)\hd\Big]\\
&&
+\hd^2(1-\frac{4}{\hs})+\hd(-8-4(\hs-5))-32\\
&=&\frac{2}{\hs}\hd^3-\frac{3\hs^4-12\hs^3+22\hs^2+2\hs+15}{2\hs^2}\hd^2\\
&&
-\frac{9\hs^3-16\hs^2-23\hs-30}{2\hs}\hd
-\frac{\hs^4-5\hs^3-\hs^2+5\hs+64}{2}
\end{eqnarray*}
(for $\hs=2$ the $-\Big(\frac{(\hs-5/2)(\hs-1)^2}{2\hs}\Big)\hd$ term
should be omitted). It is easy to compute that this implies $\hd\le
91$ for $\hs\le 6$, but $\hs=7$ we obtain only $\hd\le 153$, so taking
$r=9$ does not improve the overall bound. Taking $r=7$, we find the
overall bound $\hd\le 140$.

Generally we see from~\eqref{chi_lower_bound} that for large weights,
and hence large $\hs$, the two biggest terms in absolute value in the
cubic will be the $\hd^3$ term and a term
$-\frac{\hs^2}{4}\hd^2$. Therefore the bound on $\hd$ will be around
$|w|^3/8$.  \medskip

\noindent\emph{Example 2: weights $(1,1,1,2,6)$}

As a further example, we calculate a bound for weights
$(1,1,1,2,6)$. In this case the possible singularities are: up to
$\hd$ order~$2$ singularities (ordinary nodes) along $x_0=x_1=x_2=0$,
with $r_x=h_x=2$, and one singularity of order $6$ at $(0:0:0:0:1)$,
with $r_x=6$, $h_x=12$. At the double points, $\Delta_x^2=0$, and at
the $6$-fold point one has in fact $\Delta_x^2\ge -\frac{8}{3}$.

In this case we have $K_{\hX}^2=12K_x^2-36\hd+12\hdelta$, and
$c_1^2(\tX)=K_X^2-\sum_x \Delta_x^2\ge K_x^2-\frac{8}{3}$, so
$\theta_1=32-36\hd+12\hdelta$. We also have
\begin{eqnarray*}
c_2(\hX)&=&12c_2(\tX)+6(\hd+\hdelta)-12\sum_{x\in\cQ}(r_x-1)
-\sum_{x\in\cQ}(h_x-1)\\
&\ge&12c_2(\tX)+6(\hd+\hdelta)-12\hd-60-\hd-11\\
&=&2c_2(\tX)-71-7\hd+6\hdelta.
\end{eqnarray*}
Thus  $\theta_2=71+7\hd-6\hdelta$. Therefore
$k_0'=103$, $k_1'=-29$ and $k_2'=6$, and the quadratic is
$$
0\ge \left(1-\frac{11}{r}\right)\hd^2-(11r-274)\hd-175
$$
as long as $r\ge\hs\ge 12$ and $r^2<\hd$. Taking $r=12$, we see that
$\hd\le 699$ in this case.

We must also deal with the cases $\hs<12$ by using the cubic. For
$\hs=11$ we obtain $\hd\le 710$: as this is already bigger than $699$
it is no use looking at other choices for $r$. Smaller values of $\hs$
give smaller bounds, so the overall bound remains $\hd\le 710$.

\bibliographystyle{alpha}

\begin{thebibliography}{I-F}
\bibitem[BF]{BF} R. Braun, G. Fl{\o}ystad, \emph{A bound for the
  degree of smooth surfaces in $\PP^4$ not of general type.}
  Compos. Math. {\bf 93} (1994), 211--229.   
\bibitem[DS]{DS} W. Decker, F.-O. Schreyer, \emph{Non-general type
  surfaces in $\PP^4$: some remarks on bounds and constructions.}
  J. Symbolic Comput. {\bf 29} (2000), 545--582.
\bibitem[Do]{Do} I. Dolgachev, \emph{Weighted projective varieties.}
  In: Group actions and vector fields (Vancouver, B.C., 1981), 34--71,
  Lecture Notes in Math. {\bf 956}, Springer, Berlin 1982. 
\bibitem[EP]{EP} G. Ellingsrud, C. Peskine, \emph{Sur les surfaces
  lisses de $\PP_4$.}  Invent. Math. {\bf 95} (1989), 1--11.
\bibitem[I-F]{IF} A. Iano-Fletcher, \emph{Working with weighted
  complete intersections.} In: Explicit birational geometry of
  $3$-folds, 101--173. LMS Lecture Notes {\bf 281}, Cambridge
  University Press, Cambridge 2000.
\bibitem[Ra]{Ra} L.V. Rammea, \emph{Computations and bounds for
surfaces in weighted projective four–spaces.} Ph.D. thesis, Bath 2009.
\bibitem[Ro]{Ro} L. Roth, \emph{On the projective classification of
  surfaces.} Proc. London Math. Soc. (2) {\bf 42}, 142-170 (1936).
\end{thebibliography}

\bigskip

\noindent
L.V.~Rammea\\
Department of Mathematics and Computer Science\\
The National University of Lesotho\\
PO Roma 180\\
Lesotho\\
{\tt lv.rammea@nul.ls}\\
\smallskip

\noindent
G.K.~Sankaran\\
Department of Mathematical Sciences\\
University of Bath\\
Bath BA2 7AY\\
England\\
{\tt gks@maths.bath.ac.uk}

\end{document}